\documentclass[11pt,german,isolatin1.sty]{article}
\usepackage{amsmath,amssymb}
\begin{document}

\begin{center} {\Large\bf Do manifolds have little symmetry?}  \end{center}
\begin{center}V. Puppe\end{center}

\vspace{1cm}


\vspace{1cm}

{\bf Introduction.}\\

This note is surveying certain aspects of the following problem
stated
 by F. Raymond and R. Schultz
(s. [BH], p. 260).\\

''It is generally felt that a manifold 'chosen at random' will have very
little symmetry. Can this intuitive notion be made more precise? In
connection with this intuitive feeling, we have the following specific
question.\\
Does there exist a closed simply connected manifold, on which no
finite group acts effectively? (A weaker question,
no involution?)''\\

The problem is also listed as one of five conjectures in the
survey [AD]. But since in that reference very little information
is given about this problem, I hope that the present note can
serve as a useful complement.\\

As a general assumption we consider here closed, connected (topological)
manifolds with continuous group actions.
\vspace{1.0cm}

{\bf 1. Some remarks about non-simply connected asymmetric
manifolds.}\\

In the beginning of the 1970's several people have shown the existence
of asymmetric (i.e. not admitting any effective action of a finite group)
manifolds (see, e.g. [CRW], [B]). In fact, e.g. R. Schultz
has shown that, for dimension $\;\ge 4\;$, in any cobordism class there
are infinitely many asymmetric manifolds (s. [Sc1], [Sc2])$\;$; and several authors have given
examples of asymmetric 3-dimensional manifolds (cf. [E1]). (It is, of course, easy to see that there are no asymmetric
manifolds of dimension 1 or 2.)

All of these examples have
non-trivial fundamental group. An essential tool for many of
 the examples mentioned above is the following
result due to A. Borel, which uses the fundamental group to detect
asymmetry.
\vspace{1.0cm}

{\bf Theorem (Borel)}
If $\; M\;$ is an aspherical (i.e., $\pi_i(M) = 0$ for $i \ge 2)$
manifold, such that
\begin{itemize}
\item[(i)] $\pi_1(M)\;$ is centerless
\item[(ii)] $ Out(\pi_1(M))\;,\;$ the outer automorphism group of
$\; \pi_1(M)\;,\;$ is torsion free,
\end{itemize}
then $\; M\;$ is asymmetric (cf. [Bo],Cor.2 and [CR], Thm.3.2).\\

R. Waldm\"uller found the first example of a centerless Bieberbach
group B (i.e., a torsion free subgroup of the group of isometries
$\;{\Bbb R}^n\;,$ such that $\; M : = {\Bbb R}^n/B\;$ is a compact
aspherical manifold) with $\; Out (B) = \{1\}$ (s. [Wd]). Hence
$\;M\;$ is asymmetric by the above theorem. \vspace{1.0cm}

{\bf 2. How to use cohomology to detect asymmetry}\\

At first glance it might seem unlikely that cohomological information
could suffice to detect asymmetry. Of course, an action of a finite
group $\;G\;$ on a manifold $\;M\;$ induces an action of $\;G\;$
on the cohomology algebra $\;H^*(M)\;,$ which clearly could be trivial
without the original action being so. Hence the question is, how a non-trivial,
but cohomologically trivial action, is reflected in cohomology.\\

We first consider involutions, i.e., $\; G \cong {\Bbb Z}/2{\Bbb
Z}\;$ and we use cohomology with coefficients in $\;k = {\Bbb
F}_2$. In case of a cohomologically trivial $G-$action on $\;M\;$
the $\;E_2-$term of the Serre spectral sequence of the Borel
construction $\;M \rightarrow EG \times_G M \rightarrow BG\;$ is
isomorphic to $\;k[t] \otimes H^*(M;k)\;,\;\deg(t) = 1.$\\ The
first non-trivial higher differential is given by a derivation\\
$\;\partial : H^*(M;k) \longrightarrow H^*(M;k)\;$ of negative
degree with $\; \partial^2 = 0$. If all higher differentials
vanish, the equivariant cohomology $\;H^*_G(M;k)\;$ is isomorphic
to $\;k[t] \otimes H^*(M)\;$ as $\;k[t]-$module, but not
necessarily as $\;k[t]-$algebra. The famous Localization Theorem
for equivariant cohomology (see, e.g. [AP]) then implies that
there is a filtration on the cohomology $\;H^*(M^G;k)\;$ of the
fixed point set $\;M^G\;$ such that the associated graded algebra
is isomorphic to $\;H^*(M;k)\;$ (s. [Pu3], p. 131/132). This means
in particular that $\;H^*(M^G;k)\;$ (as a filtered algebra with
filtration $\;{\cal F}_i (H^*(M^G;k)): = \oplus^i_{j = 0}
H^j(M^G;k))\;$ is a deformation of
negative weight of the graded algebra $\;H^*(M;k)\;$ (s. [Pu2]).\\
If this deformation is trivial then $\;H^*(M;k)\;$ and
$\;H^*(M^G;k)\;$ are isomorphic as algebras yet not necessarily as
graded algebras. But if $\;H^*(M;k)\;$ has ''minimal formal
dimension'', i.e. any graded algebra (occuring as the cohomology
algebra of a manifold) , which is isomorphic to $\;H^*(M;k)\;$, as
algebra (ignoring the grading), has formal dimension bigger or
equal to that of $\;H^*(M)\;,\;$ then Smith theory implies that
the inclusion $\; M^G \longrightarrow M\;$ induces an isomorphism
of graded algebras (cf. [AHsP]). Hence in this
case $\;M^G = M\;,\;$ i.e. the action is trivial.\\

Putting all this together, we get the following result, which in a sense
is analogous to Borel's result above, but uses the cohomology
 algebra instead of the
fundamental group to detect asymmetry.
\vspace{1.0cm}

{\bf Theorem 1.} Let $\;M\;$ be a compact manifold such that
\begin{itemize}
\item[(i)$\quad$] $H^*(M;k)\;$ has no automorphism of order 2
\item[(ii)$\quad$] $H^*(M;k)\;$ has no non-trivial derivation
 of negative degree
\item[(iii)$\quad$] $H^*(M;k)\;$ has no non-trivial deformation
 of negative weight
\item[(iv)$\quad$] $H^*(M;k)\;$ has minimal formal dimension,
then $\;M\;$ does not admit any non-trivial involution.
\end{itemize}

{\bf Remark 1.} The condition (iii) in Theorem 1 can be replaced by \\
(iii)' $\; H^*(M;k)\;$ can not be given as the associated graded algebra
of a filtration of a product of Poincar\'e algebras of formal dimension
$\;< \dim \;M\;$ (see above).\\

An analogous result holds for $\;{\Bbb Z}/p{\Bbb Z}-$actions, $p\;$prime,
choosing $k = {\Bbb F}_p$.\\

Now the question is whether there exist examples fulfilling (i) - (iv)
in Theorem 1, and one might expect that they are even 'generic'.
In fact, this is true in a certain sense. If we consider 3-dimensional manifold
there ist no classification available, but we can use the
'parametrization' by their cohomology algebras (with $\;k = {\Bbb F}_2\;$
coefficients), which correspond to trilinear, symmetric forms on
$\;H^1(M;k)\;,\;$ to give a meaning to terms like 'generic' or
'chosen at random'. Or, at least to say, what is meant by
 'most ${\Bbb F}_2-$cohomology types of 3-manifolds' in Theorem 2 below.\\

Namely, if $\;\dim_k\;H^1(M;k) = m\;$ then the space of trilinear,
symmetric forms on $\;H^1(M;k)\;$ is isomorphic to the space
$\;{\cal S}^3(k^m) \cong k^{\alpha(m)}, \;\; \alpha = {m+2 \choose 3}$.
Let $\;{\cal R}(m) \subset {\cal S}^3 (k^m)\;$ be the subset of forms
that can be realized by the cohomology algebras of 3-manifolds with
$\;\dim_k H^1 (M;k) = m$. According to M. Postnikov $\;{\cal R}(m) =
{\cal R}^o(m) \cup {\cal R}^n(m)$, where $\; {\cal R}^o(m) : = \{ \mu \in
{\cal S}^3 (k^m)\;; \mu (x,x,y) + \mu (x,y,y) = 0\;$ for all
$\;x,y \in k^n\}\;,\;$ and $\;{\cal R}^n(m) : = \{\mu \in {\cal S}^3 (k^m)\;;
\exists\; x^o \in k^m\;,\; x^o \not = 0\;$ such that $\;\mu(x,x,y) + \mu
(x,y,y) = \mu (x^o,x,y)\}\;$ (s. [Po]). And let ${\cal I}(m) \subset
{\cal R}(m)\;$ be the subset of forms, which can be realized by the
cohomology algebras of 3-manifolds admitting non-trivial involutions.
By $\;|{\cal A}|\;$ we denote the number of elements of a subset $\;{\cal A} \subset
{\cal S}^3(m)$.\\
Using Theorem 1 one gets the following result.\\

{\bf Theorem 2.} Most 3-manifolds do not admit a non-trivial involution;
more precisely :

\[ \lim_{m \rightarrow \infty} \hspace{0.3cm}
\frac{| {\cal I}(m)|}{|{\cal R}(m)|} \hspace{0.3cm} = 0\;. \]
\vspace{0.4cm}

See [Pu6] for details, where in particular the connection with
binary,
self-dual codes is studied.\\
Actions of $\;{\Bbb Z}/p{\Bbb Z}\;,\; p\;$ odd prime, on 3-manifold
can be treated in a similar way. As an illustration of the method we
give a very simple proof of the following result, which was proved
independently by Su (where it is somewhat hidden in [Su], Theorem (3.9)),
J.H. Przytycki and M.V. Sokolov (s. [PS], Theorem 2.1) and A. Sikora
(s. [Si], Prop.(1.7)); cf. also Example (2.9) in [AHkP].\\

{\bf Proposition.} If a closed orientable 3-manifold $\;M\;$ admits an action
of a cyclic group $\; G \cong {\Bbb Z}/p{\Bbb Z}\;$ where $\;p\;$ is an
odd prime and the fixed point set of the action is $\;S^1\;$
then $\;H_1(M;{\Bbb F}_p) \not = {\Bbb F}_p\;$.\\

{\bf Proof.} Let us assume that $\;H^1(M; {\Bbb F}_p) \cong H_1(M;{\Bbb F}_p)
\cong {\Bbb F}_p\;$ and that $\; M^G \not = \phi\;$. We will then show
 that the Serre spectral sequence of the Borel construction collapses and
hence $\;\dim H^*(M^G;{\Bbb F}_p) = \dim H^*(M;{\Bbb F}_p) = 4\;$
by the Localization Theorem, which implies the Proposition. Since
$\; \dim H^i(M;{\Bbb F}_p) = 1\;$ for $\; i = 0,1,2,3\;$ the action must
be cohomologically trivial. It remains to show that the higher differentials
in the Serre spectral sequence vanish. Since $\; E_2 \cong H^*(BG;{\Bbb F}_p)
\otimes H^*(M;{\Bbb F}_p) \cong {\Bbb F}_p[t] \otimes \Lambda (s) \otimes
H^*(M;{\Bbb F}_p)\;$. It suffices to show that $\; H^*(M;{\Bbb F}_p)\;$
does not admit any non-trivial derivations of negative degree.\\

Let $\; 1, a_1, a_2, a_3\;$ be generators of $\;H^i(M;{\Bbb F}_p)\;$
for $\; i = 0, 1, 2, 3,\;$ such that $\; a_1 \cup a_2 = a_3\;$.
Since $\; M^G \not = \phi\;,$ 1 can not be a boundary. So any derivation
$\;\partial\;$ of negative degree vanishes on $\;a_1\;$. It follows that
$\; \partial (a_3) = \partial (a_1 \cup a_2) = (\partial a_1) \cup a_2 - a_1
\cup \partial a_2 = 0\;,\;$ if $\;\partial\;$ has degree (-1), since
$\; a^2_1 = 0\;$. Also $\; \partial a_2\;$ must vanish, since $\; 0 = \partial(a^2_2)
 = 2(\partial_2) \cup a_2\;;\;$ so $\; \partial a_2 = \lambda a_1\;$
must be zero (for $\;p\;$ odd). A derivation of degree $\;(-m)\;,\;
m \ge 2\;$ must vanish on $\; a_1\;$ and $\; a_2\;$ (since 1 is not a
boundary), and hence also on $\; a_3 = a_1 \cup a_2.\hfill_\Box$\\

Of course, one can not get simply connected 3-manifolds
without symmetry by the above approach, but the method of proof does not
 refer to the fundamental group and hence could be applied to simply
connected manifolds of higher dimension. It does not work for dimensions
4 and 5, though.\\
A. Edmonds' discussion of cyclic group actions on simply connected
4-manifolds (s. [E2], [E3]) in particular implies that there are no asymmetric
ones, and in dimension 5 the information given by the cohomology algebra
of a simply connected manifold is certainly to weak to detect
asymmetry. Hence we discuss simply connected 6-manifolds in the next
section.\\
\vspace{1.0cm}

{\bf 3. Simply connected 6-manifolds}\\

Classification theorems for certain types of simply connected
6-manifolds have been given by C.T.G. Wall ([Wa]), P.E. Jupp ([J]) and
A.V. \v{Z}ubr ([Z]). For the class $\; {\cal M}\;$ of simply connected,
6-dimensional spin-manifolds $\;M\;$ with $\;H^3(M;{\Bbb Z}) = 0\;$
the following result is contained in [Wa].\\

{\bf Theorem (Wall).} The diffeomorphism classes of elements of $\;{\cal M}
\;$ correspond bijectively to isomorphism classes of invariants.
\begin{itemize}
\item [1.$\quad$] $H\;$ free ${\Bbb Z}-$module of finite rank
(corresponding to $\;H^2(M;{\Bbb Z})\;$ for $\;M \in {\cal M})$
\item [2.$\quad$] $\mu : H \times H \times H \longrightarrow {\Bbb Z}\;$
trilinear, symmetric form (corresponding to the cup product in
$\;H^*(M;{\Bbb Z}))$
\item[3.$\quad$] $P : H \longrightarrow {\Bbb Z}\;$ linear map
(corresponding to the dual of the first Pontrjagin class)
\end{itemize}

Subject to the following conditions:

\begin{itemize}
\item[(a)$\quad$] $\mu (x,x,y) \equiv \mu (x,y,y)\quad (mod\; 2)\;\;$
for $\;\; x, y \in H$
\item [(b)$\quad$] ${\cal P}(x) \equiv 4 \mu (x,x,x)\qquad (mod\; 24)\;\;$
for $\;\; x \in H.$
\end{itemize}

Similar to Section 2 we can parametrize the elements in $\;{\cal M}\;$
by the corresponding trilinear, symmetric form in $\; {\cal S}^3({\Bbb Z}^m)
\cong {\Bbb Z}^{\alpha(m)}\;,\; \alpha(m) = {m+2 \choose 3}\;$.\\
By Wall's result this is much closer to an actual classification
 up to diffeomorphism or homeomorphism than
the parametrization in Section 2.\\
Let $\;{\cal R}(m)\;$ again denote  the set of forms, which can be realized
by the cohomology of elements in $\;{\cal M}$.\\

We define the density of a subset $\;{\cal A} \subset {\cal R}(m)\;\;$ by

\[ d_m ({\cal A}) : =\quad
 \lim \sup_{N \rightarrow \infty}\quad
\frac{| {\cal A} \cap [-N,N]^{\alpha(m)}|}
{|{\cal R}(m) \cap [-N,N]^{\alpha(m)}|} \]

Using Theorem 1 and its analogue for
 $\; G = {\Bbb Z}/p{\Bbb Z}\;,\;p\;$prime, one obtains the
 following result (s. [Pu5] for details).\\

{\bf Theorem 3.}
\begin{itemize}
\item [(a)] For $\;m \ge 6\;$ the subset of $\;{\cal R}(m)\;$
 corresponding to those manifolds, which admit a cohomologically non-trivial,
orientation preserving action of a finite group has density zero.
\item[(b)] For $\;m \ge 6\;$ the subset of $\;{\cal R}(m)\;$ corresponding to
those manifolds which admit non-trivial $\;{\Bbb Z}/p{\Bbb Z}-$actions
for infinitely many primes $\;p\;$ has density zero.
\item[(c)] For a given prime $p\;$, let $\;{\cal C}_p(m) \subset {\cal R}(m)
\;$ denote the subset corresponding to those manifolds, which admit
 a non-trivial, orientation preserving $\;{\Bbb Z}/p{\Bbb Z}-$action. Then

\[ \lim_{m \rightarrow \infty} d_m({\cal C}_p(m)) = 0.\]
\end{itemize}

Theorem 3 gives a precise meaning to the somewhat vague statement
that most manifolds in $\;{\cal M}\;$ have little symmetry.\\

{\bf Example (Iarrobino).} The following polynomial of degree 3 in
6 variables gives a trilinear, symmetric form, and hence a
Poincare duality algebra over $ {\Bbb Z} $.

$\ f(x_1,...,x_6) = 6(x_1x_4^2 - x_1^2x_4 +x_2x_4^2 + x_2x_4^2 -
 x_2^2x_5 + x_2x_5^2 + x_3^2x_4 - x_3x_4^2 + x_3^2x_6 + x_3x_6^2 +
 x_5^2x_6 + x_5x_6^2 + x_1x_2x_4 +x_1x_2x_5 + x_1x_3x_6 +x_2x_4x_6
 + x_3x_5x_6 + x_4x_5x_6 + x_4x_5x_6 + x_4^3 +x_6^3)$\\

G. Nebe has verified that this algebra has no orientation
preserving automorphisms of finite order (cf. condition (i) in
Theorem 1) and T. Iarrobino and A. Suciu have verified that it has
no deformations of negative weight modulo any prime (cf. condition
(iii) with the help of computer calculations. Conditions (ii)
 (for derivations of odd degree) and (iv) (for a non-trivial
trilinear form and $\;p\;$ odd) are easily seen to hold for all
elements in $\;{\cal R}(m)$. But clearly (iv) is not fulfilled for
$\;p = 2$. This does not matter in the case at hand if
 one assumes the $\;{\Bbb Z}/2{\Bbb Z}-$action to be orientation preserving
(because then the fixed point set must have even codimension),
but it shows that one can not exclude the possibility of orientation
reversing involutions by applying Theorem 1. So one gets the following
result.
\vspace{1.0cm}

{\bf Theorem 4.} There exist simply connected manifolds on which no
finite group can act effectively and orientation preserving.\\
(See [Pu5]).\\

For orientation reversing involutions one can imitate the
arguments leading to Theorem 3 up to a certain point. This is
sketched in [Pu5], Remark 1.3. But since some of the consequences
are prerequisites for Kreck's result (s.[K]),
 we give some more details here. \\

Any algebra $H^*(M, \mathbb{Z}), M \in \mathcal{M}$, admits an
orientation reversing involution, given by the identity on $H^0$
and $H^4$, and by multiplication with $(-1)$ on $H^2$ and $H^6$,
in contrast to the fact that most $H^*(M; \mathbb{Z})$ do not
admit any non-trivial orientation preserving involution. Since the
composition of an arbitrary orientation reversing involution with
the particular one given above is an orientation preserving
involution, it follows that most $H^*(M; \mathbb{Z})$ only admit
the one orientation reversing involution given
above.\\

The following results imply for most $M \in \mathcal{M}$, if
$\tau: M \to M$ is a non-trivial involution, then $(M, \tau)$ is a
conjugation space in the sense of [HHP], in particular $H^*(M;
\mathbb{F}_2) \cong H^*(M^{\tau}; \mathbb{F}_2)$ by a degree
halving isomorphism of
algebras, $M^{\tau}$= fixed point set. (cf. [Pu5], Remark 1.3)\\

{\bf Theorem 5.} Let $M \in \mathcal{M}$, and let $H^*(M;
\mathbb{F}_2)$ be generated by $H^2(M; \mathbb{F}_2)$ as an
algebra. Assume that $\tau: M \to M$ is an involution with
non-empty fixed point set $M^{\tau}$, which acts on $H^2(M;
\mathbb{Z})$ by multiplication with $(-1)$. Then $(M,\tau)$ is a
conjugation space.\\

{\bf Proof.} We first consider the Serre spectral sequence of the
Borel construction of the $C$-space $M$, $C: = \{id, \tau\} \cong
\mathbb{Z}/{2\mathbb{Z}}$ with coefficients in $\mathbb{F}_2$.
Since $\tau$ acts trivially on $H^*(M; \mathbb{F}_2)$ (which
follows immediately from our assumption),
 we have
 $$
 E_2^{*,*} \cong H^* (BC; \mathbb{F}_2) \otimes H^* (M; \mathbb{F}_2).
 $$
The differentials in the spectral sequence correspond to
derivations of negative degree on $H^*(M; \mathbb{F}_2)$. Since
$\mathcal{M}^{\tau}$ is assumed to be non-empty, all these
derivations vanish on $H^2(M; \mathbb{F}_2)$. (Otherwise $H^*_C(M;
\mathbb{F}_2)$ would be annihilated by a power of $u \in
\mathbb{F}_2 [u] \cong H^* (BC; \mathbb{F}_2),  deg(u)=1$, and
hence $M^{\tau}$ would have to be empty by the Localization
Theorem.) Since $H^*(M; \mathbb{F}_2)$ is generated by $H^2(M;
\mathbb{F}_2)$, all differentials in the spectral sequence vanish,
and the inclusion $M^{\tau} \to M$ induces an injection $H^*_C(M;
\mathbb{F}_2) \to H^*_C(M^{\tau}; \mathbb{F}_2) =
H^* (BC; \mathbb{F}_2) \otimes H^*(M^{\tau}; \mathbb{F}_2)$ (see, e.g. [AP], Prop. (1.3.14)).\\

We next try to calculate the analogous map for $\mathbb{Z}$
coefficients. The $E_2$-term of the Serre spectral sequence for
the Borel construction of $M$ is given by $E^{*,*}_2 \cong H^*(BC;
H^*(M; \mathbb{Z}))$. The map induced by coefficient change,
$\mathbb{Z} \to \mathbb{F}_2$, on the $E_2$-terms is injective for
all $E^{i;j}_2$ with $i>0$. So the spectral sequence with
$\mathbb{Z}$ coefficients also collapses at the $E_2$-level.
Since, up to periodicity (i.e. multiplication by $t \in \mathbb{Z}
[t] /(2t) \cong H^*(BC; \mathbb{Z}),  deg (t)=2)$, the only
non-zero terms in the $E_2$-term are $E^{4,0}_2, E^{1,2}_2$ and
$E^{1,6}_2$, it is easy to see that there are no extension
problems for $E_2= E_x$. This is shown in a more general context
by M. Olbermann (s.[O1]).\

Hence, as $H^*(BC; \mathbb{Z}$-module, $H^*_C(M; \mathbb{Z})$ is
isomorphic to $H^*(BC; \mathbb{Z}) \otimes (H^0 \oplus H^4) \oplus
m \otimes (H^2[-1] \oplus H^6 [-1])$, where $m: = ker (H^*(BC;
\mathbb{Z}), \to \mathbb{Z}), t \mapsto 0,H^i:= H^i (M;
\mathbb{Z})$, and $''[-1]''$ indicate a degree shift by $(-1)$,
i.e. we identify $t \otimes H^2(M; \mathbb{Z}) [-1]$ with $
E^{1,2}_2 = H^1 (BC; H^2 (M; \mathbb{Z}))$, etc. Here, the
assumption is used that $\tau$ acts on $H^2$ (and $H^6$) by
multiplication with $(-1)$.\

Next we calculate the equivariant cohomology of the fixed point
set $M^{\tau}$ with $\mathbb{Z}$ coefficients. This is not
completely obvious since $H^*(M^{\tau}; \mathbb{Z})$ may have
$\mathbb{Z}$-torsion. By the Universal-Coefficients-Theorem we get
an inclusion of algebras
$$
H^*_C(M^{\tau}; \mathbb{Z}) \otimes \mathbb{F}_2 \longrightarrow
H^*_C(M^{\tau}, \mathbb{F}_2)=H^*(BC; \mathbb{F}_2) \otimes
H^*(M^{\tau}; \mathbb{F}_2).
$$
The image of this map is contained in the kernel of the Bockstein
operator $\beta: H^*_C(M^{\tau}; \mathbb{F}_2) \to M^{*+1}_c
(M^{\tau}; \mathbb{F}_2)$, more precisely: The intersection of
this image with the kernel of the restriction to the fibre
$H^*_C(M^{\tau}; \mathbb{F}_2) \to H^*(M^{\tau}; \mathbb{F}_2)$ is
the subalgabra of $H^*(BC; \mathbb{F}_2) \otimes H^* (M^{\tau};
\mathbb{F}_2)$ given  by all elements of the form $u^{2k} \otimes
x + u^{2k-1} \otimes \beta (x)$, where $x \in H^*(M^{\tau};
\mathbb{F}_2)$ and $u^{2k}, u^{2k-1} \in H^*(BC; \mathbb{F}_2)
\cong \mathbb{F}_2 [u]; k>0$. (Recall that $\beta (u^{2k}) =0$ and
$\beta (u^{2k-1}) = u^{2k}.)$ \

We consider the following commutative diagram
$$
\begin{array}{ccc }
   H^*_C (M; \mathbb{Z})   &   \stackrel{r}{\longrightarrow} &  H^*_C (M^{\tau}; \mathbb{Z})  \\
    \downarrow &  & \downarrow \\
    H^*_C(M; \mathbb{F}_2) & \stackrel{\bar{r}}{\longrightarrow} &    H^*_C(M^{\tau}; \mathbb{F}_2)
\end{array}
$$
where the horizontal maps are induced by the inclusion $M^{\tau}
\to M$, and the vertical maps by coefficient change $\mathbb{Z}
\to \mathbb{F}_2$.\

An element of the form $t \otimes a \in t \otimes H^2(M;
\mathbb{Z}) [-1] \subset H^*_C(M; \mathbb{Z})$ is mapped to $u
\otimes \bar{a} \in H^*_C (M; \mathbb{F}_2) \cong H^*(BC;
\mathbb{F}_2) \stackrel{\sim} {\otimes} H^*(M; \mathbb{F}_2)$
under coefficient change, where $\bar{a} $ is the image of $a$ in
$H^2(M; \mathbb{F}_2)$. Since the restriction of $t \otimes a$ to
the fibre, $H^*(M; \mathbb{Z})$, vanishes, the image of the
element $r(t \otimes a)$ in $H^*_C (M^{\tau}; \mathbb{F}_2)$ is of
the form $u^2 \otimes x + u \otimes \beta(x)$, with $x \in
H^1(M^{\tau}; \mathbb{F}_2)$. It follows that $\bar{r}(1 \otimes
\bar{a}) = u \otimes x + 1 \otimes \beta(x)$. Since $H^2 (M;
\mathbb{F}_2)$ generates $H^*(M; \mathbb{F}_2); 1 \otimes H^2(M;
\mathbb{F}_2)$ generates $H^*_C (M; \mathbb{F}_2)  \cong H^* (BC;
\mathbb{F}_2)  \stackrel{\sim}{\otimes} H^* (M; \mathbb{F}_2)$ as
$H^*(BC, \mathbb{F}_2)-$algebra, even though the multiplication
might - a priori - be twisted (which is indicated by the tilde
sign). One therefore gets that $H^*(M^{\tau}; \mathbb{F}_2)$ is
generated, as $\mathbb{F}_2-$ algebra, by those $x \in
H^1(M^{\tau}, \mathbb{F}_2)$ which occur in $\bar{r}(1 \otimes
\bar{a})$ for $\bar{a} \in H^2(M^{\tau}, \mathbb{F}_2)$ (cf.
[Pu1]). This implies that the composition
$$
H^*_C (M; \mathbb{F}_2) \to H^*_C(M^{\tau}; \mathbb{F}_2)=H^*(BC;
\mathbb{F}_2) \otimes M^*(M^{\tau}; \mathbb{F}_2)
$$
$$
\rightarrow H^*(BC; \mathbb{F}_2) \otimes H^*(M^{\tau};
\mathbb{F}_2)/ \oplus _{i<j} H^i (BC; \mathbb{F}_2) \otimes H^j
(M^{\tau}; \mathbb{F}_2)
$$
is an isomorphism of $\mathbb{F}_2$-vector spaces. This property
is the dual equivalent of a characterization of conjugation spaces
due to M. Olbermann (s. [O2]).
Hence $(M, \tau)$ is a conjugation space. $\hfill_\Box$\\

{\bf Remark 2.} The proof of the above theorem carries over from
$M \in \mathcal{M}$ to manifolds $M$ of arbitrary (even) dimension
with
$H^{odd}(M;\mathbb{Z}) = 0$. \\

The following Lemma implies that for most $M \in \mathcal{M}$ the
algebras $M^*(M; \mathbb{F}_2)$ do not admit any non-trivial
derivation of negative
degree (cf. [Pu6], Remark 5).\\

{\bf Lemma 1.}  Let $A^*= A^0 \oplus A^2 \oplus A^4 \oplus A^6$ be
a graded connected Poincar\`e duality algebra over $\mathbb{F}_2$
which is generated by $A^2$. The algebra $A^*$ admits a
non-trivial derivation of negative degree if and only if there
exists a subspace  $K \subset A^2$ of codimension $1$ with the
following properties:
\begin{itemize}
  \item[(i)]$K \times K \stackrel{\mu}{\to} K^{\bot} \subset A^4$,
  where $K^{\bot}$ denotes the orthogonal complement with
  respect to the Poincar\`e pairing
  \item[(ii)] $\exists a \in A^2 \backslash K$ such that
\begin{itemize}
  \item[(1)] the map $K \stackrel{\mu_a}{\to} A^4 \to A^4 \backslash K^{\bot}$ is
  an isomorphism where $\mu_a: K \to A^4$ denotes
  the multiplication by $a$.
  \item[(2)] $aa \in K^{\bot}$.
\end{itemize}
\end{itemize}

{\bf Proof.} Assume that  $\partial: A^* \to A^*$ is a non-trivial
derivation of negative degree. Put $K:= ker \partial|_{A^2}$.
Condition (i) means that $k_1 k_2 k_3=0$ for all $k_i \in K,
i=1,2,3$. Clearly $\partial (k_1 k_2 k_3) = 0$, since $\partial
(k_i)=0$ for $i=1,2,3$. But because of Poincar\`e duality (and $1
\in \partial (A^2)$), $\partial|_{A^6}$ is injective.
This shows (i). \\

For any $k \in K$ one has $\partial (aak) = \partial (aa)k +
aa\partial k= 0$. Hence $aak=0$, and $aa \in K^{\bot}$. If $k
\not= 0$ then, by Poincar\`e duality there exists a $k_1 \in k$
with $aak \not= 0$; and therefore for any $x \in A^2, x \not= 0$
there exists a $k \in K$ such that $xk \not= 0$. If $b \in
K^{\bot}$ and $k \in K$, then $0 = \partial (bk) = \partial (b)k$.
Hence $\partial(b)= 0$. So, $\partial : A^4 \to A^2$ factors
through $A^4/K^{\bot}$ and for the composition $K
\stackrel{\mu_a}{\to} A^4 \to A^4/K^{\bot} \to K$ is the identity
since $\partial(ak) = \partial(a)k = k$. Hence also (ii) holds. \\

Assume now that (i) and (ii) are fulfilled. Define
$\partial:A^*\to A^*$ by $\partial^2:A^2 \to A^2/K \cong
\mathbb{F}_2 = A^0$, $\partial^4:A^4 \to A^4/K^{\bot}
\stackrel{{\mu_a}^{-1}}{\rightarrow} K \subset A^2$, $
\partial^6:A^6 \stackrel{\cong} \to  K^{\bot} \subset A^4$. Clearly
$\partial^*
\partial^* = 0 $. We still have to check the derivation property.
For $k_1,k_2 \in K$ we have $ 0 = \partial (k_1 k_2) = (\partial
k_1)k_2 + k_1(\partial k_2)$; for $a$ and $k \in K$ one gets
$\partial (ak) = (\partial a)k + a(\partial k) = k$ by definition
of $\partial^4$, $\partial (aa) = (\partial a)a + a(\partial a) =
0$, since $aa \in K^{\bot}$. Using the decomposition $A^2 = K
\oplus <a>$ and $A^4 = K^{\bot}\oplus <a>^{\bot}$, where $<a>$ is
the ${\mathbb F}_2$-vector space generated by $a\in A^2 \backslash
K$, it is straight forward to check for all products that the
derivation property holds.$\hfill_\Box$\\

{\bf Remark 3.} In the presents of condition (i) of the above
lemma conditions(ii), (1) and (2), are equivalent to conditions
(1) and (2) for all $a \in A^2\backslash K$.\\

For the example above a straight forward calculation gives the
following structure constants for the trilinear form mod 2, $\mu
:{\mathbb F}^6_2 \times {\mathbb F}^6_2 \times {\mathbb F}^6_2 \to
{\mathbb F}^6_2$ ,$\mu _{ijk}:= \mu (e_i,e_j,e_k), i,j,k =
1,...,6$ , corresponding to the cubic polynomial with respect to
the coordinate system $x_1, ..., x_6$: $\mu _{124} = \mu _{125} =
\mu _{136} = \mu _{246} = \mu _{356} = \mu _{456} = 1$; and all
other $\mu _{ijk}$ (not obtained by permutation of indices from
the former) vanish.\\

It is easy to check in this example that, mod 2, $A^2$ generates
$A^*$, and also that the conditions in Lemma 1 for the existence
of
a non-trivial derivation of negative degree are not fulfilled.\\


{\bf Corollary.} If $(M,\tau)$ is a manifold with an orientation
reversing involution, $M \in \mathcal{M}$, and $H^*(M;
\mathbb{Z})$ given
by the polynomial in the above example, then $(M, \tau)$ is a conjugation space.\\

Recently M. Kreck (s. [K]), using completely different methods,
has shown that in the situation of the above corollary the first
Pontrjagin class of $M$ has to fulfill an additional condition, if
the action is differentiable, or - at least - locally linear. He
obtained the following result, which completely
answers the question, stated in the introduction in the smooth and locally linear category.\\

{\bf Theorem (Kreck).} There are infinitely many closed asymmetric simply connected smooth 6-manifolds.\\

The above Theorem (Wall) is only a part of Wall's classification result in
that he considers the bigger class $\;{\cal N}\;$ of simply connected,
6-dimensional spin-manifolds with free integral cohomology, so
$\; H^3(M;{\Bbb Z})\;$ is a free module of even rank (because of
Poincar\'e duality) for $\; M \in {\cal N}\;.$ Wall reduces the
classification of $\;{\cal N}\;$ to that of $\;M\;$ by showing that
a manifold $\; M \in {\cal N}\;$ is diffeomorphic to $\; M^\prime \sharp
S^3 \times S^3 \cdots \sharp S^3 \times S^3\;$ with $\; M^\prime \in
{\cal M}\;.$ So the only additional invariant needed is the rank of
$\; H^3(M;{\Bbb Z})\;.$\\

It is clear that no manifold $\; M \in {\cal N}\;$ with $\; H^3(M;{\Bbb Z})
\not = 0\;$ fulfills assumption (i) of Theorem 1, since then $\; H^*(M;{\Bbb F}_2)
\;$ admits non-trivial involutions (as graded algebra).
Hence for the following we restrict to cohomologically trivial actions,
i.e. the induced action on $\; H^*(M;{\Bbb Z})\;$ is assumed to be trivial.\\
Considering the Serre spectral sequence of the Borel construction
with coefficients in $\;{\Bbb Z}\;$ (cf. [Pu5], Prop. 1) the
assumption (ii) in
 Theorem 1 is fulfilled for $\; M \in {\cal M}\;$ already for degree
reasons. But for elements in $\; {\cal N}\;$ one needs an extra argument,
namely the following simple lemma.\\

{\bf Lemma 2.} Let $\; A^*\;$ be a Poincar\'e duality algebra over
$\; {\Bbb F}_p\;,\;p\;$ prime, of formal dimension 6, and let $\;
A^1 = A^5 = 0\;.\;$ Assume that the even dimensional part $\;
A^{ev}\;$ is generated by $\;A^2\;$ (as an algebra with unit),
then $\;A^*\;$ does not admit a
non-trivial derivation of negative, odd degree.\\

{\bf Proof.} Let $ A^0 = <1>\;,\; A^2 = <a_i\;,\; i \in I>\;,\;A^3 =
<b_j\;,\;\bar{b}_j\;;\;j \in J>,\\  A^4 = <c_i\;,\;i \in I>\;,\;
A^6 = <d>\;,\;$ where the $\; \bar{b}_j\;,\;c_i\;,\;d\;$ form the dual basis of
$\;b_j\;,\;a_i\;,\;1\;$ with respect to the Poincar\'e duality pairing,
 and $<\;\;>$ denotes the vector space generated by the indicated basis.
Assume that $\;\partial : A^* \longrightarrow A^*\;$is a derivation of
degree (-1). Then $\;\partial\;$ vanishes on $\;A^0\;,\;A^2\;$ and
$\;A^6\;$ for degree reasons, and on
$\;A^4\;$ since $\;A^{ev}\;$ is generated by $\;A^2\;$.\\
Let $\;b\;$ be a non-zero element in $\;A^3\;$ and let $\;a : = \partial b\;$.
Assume that $\;a \not = 0\;$. If $\;c\;$ is the dual of $\;a\;$, so
$\;ac = d\;,\;$ then $\;\partial(bc) = (\partial b)c - b\;\partial c = ac = d\;;\;$
but $\;bc = 0\;.$ Hence we get a contradiction. So $\;a = \partial b = 0\;$ for any $\;
b \in A^3\;.$ Therefore $\;\partial\;$ must be trivial.\\
 The argument for
derivations of lower (negative, odd) degree
 is completely analogous.\hfill$_\Box$\\

{\bf Remark 4.} It is obvious that one can generalize Lemma 1 to
Poincar\'e
 duality algebras $\;A^*\;$ of formal dimension $\;2m\;$ with $\;A^1 = 0\;$
and $\;A^{ev}\;$ generated by $\;A^2\;$.\\

Parametrizing $\; {\cal N}\;$ by integral cohomology type one gets the
following generalization of Theorem 3.\\

{\bf Theorem 6.} Most integral cohomology types in $\;{\cal N}\;$
do not admit non-trivial but cohomologically trivial
 $\;{\Bbb Z}/p{\Bbb Z}-$actions.\\

Here 'most' can be given a precise meaning similar to Theorem 3.\\

For the proof one uses the integral version of Theorem 1 for
$\;{\Bbb Z}/p{\Bbb Z}-$actions. By assumption the considered actions
are cohomologically trivial and\\
$\; H^*(M;{\Bbb Z})\;$ is free.
So the $\; E_2-$term of the Serre spectral sequence of the Borel
construction with integral coefficients is given by

$E_2 \cong H^*(BG;H^*(M;{\Bbb Z}) \cong H^*(BG;{\Bbb Z}) \otimes
H^*(M;{\Bbb Z}) \cong {\Bbb Z}[t]/(pt) \otimes H^*(M;{\Bbb Z})\;,$\\
where $\deg(t) = 2.$

The first non-trivial boundary in the spectral sequence would give a
non-trivial derivation of odd, negative degree on $\; H^*(M;{\Bbb Z})
\otimes {\Bbb F}_p \cong H^*(M;{\Bbb F}_p)\;$. Hence, by Lemma 1,
the spectral sequence collapses if $\; H^{ev}(M;{\Bbb F}_p)\;$ is
generated by $\; H^2(M;{\Bbb F}_p)\;,$ which is the case if the
 trilinear form,
given by the cup product, is non-degenerate. This holds for most
cohomology types. Since $\;t\;$ above has degree 2,
$\;\;H^{ev}(M^G;{\Bbb F}_p)\;$ is a deformation of negative weight
of $\;H^{ev}(M;{\Bbb F}_p)\;\;$ (cf. [Pu5], Prop. 1 and 4). In
most cases there are no such non-trivial deformations. So
$\;H^{ev}(M^G;{\Bbb F}_p)\;$ is isomorphic to $\;H^{ev}(M;{\Bbb
F}_p)\;$ as filtered algebras. In case of a non-degenerate
trilinear form the cup length of $\;H^{ev}(M;{\Bbb F}_p)\;$ (and
hence of $\;H^{ev}(M^G;{\Bbb F}_p)) \;$ is 3. So the dimension of
$\;M^G\;$ must be 6, and hence $\;M^G = M\;;\;$ i.e.
 the action is trivial.\hfill$_\Box$\\

{\bf Remark 5.} If $\;p\;$ is large compared to the size of
$\;H^*(M;{\Bbb Z}) \;$ (more precisely: $\;\;p>\; rkH^i (M;{\Bbb
Z}) + 1\;,\;$ for all $\;i\;)\;,\;$ then an action of $\;{\Bbb
Z}/p{\Bbb Z}\;$ on $\;M\;$ must be cohomologically trivial by
elementary representation theory. So Theorem 5 holds for 'large
$\;p\;$' without the restriction 'cohomologically trivial'. But
recently M.Olbermann has improved Theorem 6 considering not only
cohomologically trivial actions, but assuming that $\ rk
H^3(M;{\Bbb Z}) = o((rkH^2(M;{\Bbb Z})) ^{3/2}) $ (s.
[O1])\\

The classification of simply connected 6-manifolds, without the
assumptions 'spin' and $\;'H^*(-;{\Bbb Z})\;$ free over $\;{\Bbb
Z'}\;$ involves further invariants (s. [Z]). But for a given
manifold $\;M\;$ we can kill the torsion in $\;H^*(M;{\Bbb Z})\;$
localizing $\;{\Bbb Z}\;$ by inverting those primes which occur in
the torsion. The above arguments can then be applied to $\;{\Bbb
Z}/p{\Bbb Z}-$action, where $\;p\;$ does not belong to the
(finitely many) inverted primes. Parametrizing the class of simply
connected manifolds by their rational cohomology algebras (which
is, of course, far from a classification up
to homeomorphism or diffeomorphism) one gets the following result.\\

{\bf Theorem 7.} Most rational cohomology types of simply
connected 6-manifolds do admit non-trivial $\;{\Bbb Z}/p{\Bbb
Z}-$action for at most finitely many primes.\\

This generalizes Theorem 2 in [Pu4].\\

{\bf Remark 6.} It is to be expected that similar results hold for
higher (even) dimensions. But in particular the discussion of
condition (iii) or (iv) in Theorem 1 gets more and more involved
with increasing dimension. Certain results for $\;S^1-$actions on
8-manifolds in this direction are contained in [I].

\vspace{1cm}


{\footnotesize {\it {Volker Puppe\\
Fachbereich Mathematik und Statistik\\
Universt\"at Konstanz\\
D-78457 Konstanz\\
Germany\\}
 E-mail: volker.puppe@uni-konstanz.de}}


\begin{thebibliography}{10000}
\bibitem[AD] {AD} Adem,A.,Davis,J.F.: {\it Topics in transformation
groups}, Handbook of gemetric topology, 1-54, North-Holland,
Amsterdam (2002).
\bibitem[AP]{AP} Allday, C., Puppe, V.: {\it Cohomological Methods in
Transformation Groups}, Cambridge University Press (1993).
\bibitem[AHsP]{AHsP} Allday, C., Hauschild, V., Puppe, V.: {\it A non-fixed
point theorem for Hamiltonian Lie group actions}, Trans. Amer.
Math. Soc. 354, 2971 -2982 (2002).
\bibitem[AHkP]{AHkP} Allday, C., Hanke, B., Puppe, V.: {\it Poincar\'e
duality in P.A. Smith theory}, Proc. Amer. Math. Soc. 131, 3275
-3283 (2003). .
\bibitem[Bl]{Bl} Bloomberg, E.M.: {\it Manifolds with no periodic
homeomorphism}, Trans. Amer. Math. Soc. 202, 67-78 (1975).
\bibitem[Bo] {Bo} Borel, A.: {\it On periodic maps of certain $
K(\pi ,1)$}, Collected papers III, 57-60, Springer (1983).
\bibitem[BH]{BH} Browder, W., Hsiang, W.C.: {\it Some problems on
homotopy theory, manifolds and transformation groups}, Proc. Symposia
Pure Math. 32, 251-267 (1978).
\bibitem[CRW]{CRW} Conner, P.E., Raymond, F., Weinberger, P.:
{\it Manifolds with no periodic maps}, In: Proc. Second Conference Compact
Transformation Groups, Part II, Springer Lect. Notes 299, 81-108 (1972).
\bibitem[E1]{E1} Edmonds, A.: {\it Transformation groups and low-dimensional
manifolds, Group Actions on Manifolds}, Contemporary Mathematics 36, 339-
366 (1985).
\bibitem[E2]{E2} Edmonds, A.: {\it Construction of group actions on
four-manifolds}, Trans. Amer. Math. Soc. 299, 155-170 (1987).
\bibitem[E3]{E3} Edmonds, A.: {\it Aspects of group actions on four-manifolds},
Topology Appl. 31, 109-124 (1989).
\bibitem[HHP]{HHP} Hausmann, J.-C., Holm, T., Puppe, V.: {\it
Conjugation spaces}, Algebraic and Geometric Topology 5, 923-964
(2005).
\bibitem[I]{I} Iniotakis, J.M.: {\it Mannigfaltigkeiten mit wenig Symmetrie},
Diplomarbeit, Konstanz (1999).
\bibitem[J]{J} Jupp, P.E.: {\it Classification of certain 6-manifolds},
Proc. Camb. Phil. Soc. 73, 293-300 (1973).
\bibitem[K]{K} Kreck,M.: {\it Simply connected asymmetric manifolds},
preprint (2006).
\bibitem[KP]{KP} Kreck,M., Puppe, V.:{\it Involutions on
3-manifolds and self-dual, binary codes}, in preperation.
\bibitem[O1]{O1} Olbermann, M.:{\it 'Most' manifolds in ${\cal
N}$ do not admit orientation preserving $\;{\Bbb Z}/p{\Bbb Z}-$
actions}, preprint (2005).
\bibitem[O2]{O2} Olbermann, M.: {\it A definition of conjugation
spaces without a conjugation equation}, preprint (2005).
\bibitem[Po]{Po} Postnikov, M.M.: {\it Construction of intersection rings of
3-dimensional manifolds (Russian)}, Dokl. Akad. Nank SSSR 61, 795-797 (1948).
\bibitem[PS]{PS} Przytycki, J.H., Sokolov, M.: {\it Surgeries on
periodic links and homology of periodic 3-manifolds}, Math. Proc. Camb.
Phil. Soc. 131, 295-307 (2001).
\bibitem[Pu1]{Pu1} Puppe, V.: {\it On a conjecture of Bredon},
manuscripta math. 12,11-16 (1974).
\bibitem[Pu2]{Pu2} Puppe, V.: {\it Cohomology of fixed point sets and deformation
of algebras}, manuscripta math. 23, 343-354 (1978).
\bibitem[Pu3]{Pu3} Puppe, V.: {\it Deformations of algebras and cohomology
of fixed point sets}, manuscripta math. 30, 119-136 (1979).
\bibitem[Pu4]{Pu4} Puppe, V.: {\it Simply connected manifolds without
$S^1-$symmetry}, In: Topology Conference Göttingen 1987, Springer Lect. Notes
1361, 261-268 (1988).
\bibitem[Pu5]{Pu5} Puppe, V.: {\it Simply connected 6-dimensional manifolds
with little symmetry and algebras with small tangent space},
Prospects in Topology, Annals of Math. Studies 138, 283-302 (1995).
\bibitem[Pu6]{Pu6} Puppe, V.: {\it Group actions and codes}, Canad. J.
Math. 53, 212-224 (2001).
\bibitem[Sc1]{Sc1} Schultz, R.: {\it Group actions on hypertoral manifolds I},
Proc. Topology Symposium Siegen 1979, Springer Lect. Notes 788, 364-377
(1980).
\bibitem[Sc2]{Sc2} Schultz, R.: {\it Group actions on hypertoral manifolds II},
J. Reine Angew. Math. 325, 75-86 (1981).
\bibitem[Si]{Si} Sikora, A.: {\it Torus and ${\Bbb Z}_p-$actions on
manifolds}, Topology 43, 725 -748 (2004).
\bibitem[Su]{Su} Su, J.C.: {\it Periodic transformations on the product of
two spheres}, Trans. Amer. Math. Soc. 106, 305-380 (1963).
\bibitem[Wd]{Wd} Waldm\"uller, R.: {\it Eine flache Mannigfaltigkeit
ohne Symmetrien}, Diplomarbeit, Aachen (2002).
\bibitem[Wa]{Wa} Wall, C.T.C.: {\it Classification problems in differential
topology. V. On certain 6-manifolds}, Invent. math. 1, 355-374 (1966).
\bibitem[Z]{Z} \v{Z}ubr, A.V.: {\it Classification of simply connected
topological 6-manifolds}, Topology and Geometry-Rohlin Seminar,
Springer Lect. Notes 1346, 325-339 (1988).
\end{thebibliography}
\end{document}